
\mag=\magstep1
\documentstyle{amsppt}

\topmatter
\title Semistable Minimal Models of Threefolds in Positive or
Mixed Characteristic
\endtitle
\rightheadtext {SEMISTABLE MINIMAL MODELS}
\author Yujiro Kawamata \endauthor
\address Department of Mathematical Sciences, University of Tokyo, Hongo,
Bunkyo, Tokyo, 113, Japan \endaddress
\email kawamata\@math.s.u-tokyo.ac.jp\endemail

\abstract
We extend the minimal model theorem to the
3-dimensional schemes which are
projective and have semistable reduction
over the spectrum of a Dedekind ring.
\endabstract
\endtopmatter

\document

\head 0. Introduction \endhead

The classification theory of algebraic surfaces is characteristic free (cf.
[Mu]), and it also holds for log surfaces, i.e., pairs consisting of
surfaces and divisors (cf. [K1], [F], [TM]).
But the minimal model theory of algebraic 3-folds is only completed
over the field of characteristic 0,
because the vanishing theorem of Kodaira type is
necessary in the course of the proof (cf. [KMM], see also [Ko]).
The purpose of this paper is to extend the minimal model theory
to the positive or mixed characteristic case
for 3-dimensional schemes which have semistable reduction
over the spectrum of a Dedekind ring.
We refer to [KMM] for the terminology of the minimal model theory.

Because our 3-fold is fibered by a family of surfaces,
we can prove the Cone Theorem (Theorem 1.3) easily
by using the results on log surfaces.
It says that there exists an extremal ray if the given 3-fold is not
minimal relatively over the base.
Since the vanishing theorem of Kodaira type is false in positive
characteristic, we cannot prove the Base Point Free Theorem in general.
But using certain weaker vanishing theorems (Lemmas 2.1 and 2.2),
we prove the Contraction Theorem (Theorem 2.3) :
there exists a contraction morphism associated to an extremal ray.

There are several things about which we should be careful in the positive
or mixed characteristic case.
Because we do not have the Grauert-Riemenschneider vanishing theorem,
we do not know whether terminal singularities
are always Cohen Macauley.  So we put a technical assumption (6)
at (1.1) instead.
This condition is automatically satisfied in the case of characteristic 0,
but we do not know whether it is true or not in the general case.
We shall prove that this assumption is preserved under
the birational transformations in the process
toward the minimal models, i.e., the divisorial
contractions and flips (Theorems 5.3 and 5.5).

In order to prove the Flip Theorem in \S5,
we classify the {\it terminal singularities} which appear
inevitably in the course of the minimal model program in \S\S3-4.
Here we may assume that
the base ring is a complete discrete valuation ring whose residue field
is algebraically closed.
The canonical cover is an inseparable morphism if its degree is divisible
by the characteristic,
and the global canonical covering space may have non-normal points,
e.g., for a quasi-hyperelliptic surface (cf. [BM]).
But we can prove similar statements
for a local canonical cover as in the case of characteristic 0
by using a general position argument (\S3).
Then the singularities of the canonical covers are
classified by using the assumption (1.1) (\S4).
Here we can use the method of the toric geometry (cf. [KKMS]),
because it concerns the set of monomials,
and independent of the characteristic.

The termination of the sequence of flips (Theorem 5.2) is proved by
looking at the minimal discrepancy coefficients of exceptional
divisors (cf. [K4]).
Then by induction on the maximal index of the relative canonical divisor,
we prove the Flip Theorem (Theorem 5.5) :
the flip exists for a small contraction.
This proof is new and simple
compared with the previously known proofs when it is applied to
the case of characteristic 0 (cf. [T], [K2], [S], [Mo2]).
But the additional condition (6) of (1.1) made our proof longer.

As in [KMM], the combination of the Cone, the Contraction and the Flip
Theorems yields the main theorem of this paper,
the Minimal Model Theorem (Theorem 5.7) :
the given morphism has a birational model which is minimal or has a
Mori fiber space structure.

The author would like to thank Professor S. Tsunoda for explaining
[TM2] and the unpublished version of [T].

\head 1. Cone theorem
\endhead

Let $\Delta = \text{Spec } A$ be the spectrum of a Dedekind ring,
and let $f: X \to \Delta $ be a flat quasi-projective morphism
from a normal scheme $X$ of dimension 3,
or a morphism obtained by a completion from such a morphism.
In this paper, $f$ is said to be {\it semistable} if the following
conditions are satisfied:
\roster
\item"(1)"  the generic fiber  $X_{\eta}$  is smooth,

\item"(2)"  $X$ is regular, and the fibers $X_s$ for the closed points
$s \in \Delta$ are geometrically reduced normal crossing divisors,

\item"(3)"  the irreducible components of any closed fiber
$X_s$ are geometrically irreducible and smooth.
\endroster

The purpose of this paper is to construct a minimal model
for the given semi\-stable and projective morphism $f: X \to \Delta$.

The {\it relative canonical sheaf}
$\omega _{X/\Delta} = H^{- 2}(\omega _{X/\Delta}^{\bullet})$
is a coherent reflexive sheaf of rank 1.
The corresponding Weil divisor $K_{X/\Delta}$ is called the
{\it relative canonical divisor}.

The construction of a minimal model is inductive,
and we shall deal with the
category of morphisms given below.  Note that a semistable morphism is
in this category as well as its minimal model, but the latter is not
necessarily semistable.

\definition{1.1. Assumption}
Throughout this paper, we assume the following conditions:
\roster
\item"(1)"  the generic fiber  $X_{\eta}$  is smooth,

\item"(2)"  the fibers $X_s$ for the closed points $s \in \Delta$ are
geometrically reduced and satisfy the condition $(S_2)$,

\item"(3)"  the irreducible components $S_i$ of any closed fiber
$X_s$ are geometrically irreducible and geometrically normal
$\Bbb Q$-Cartier divisors on $X$,

\item"(4)"  there exists a finite set $\Sigma$ of closed points on $X$
such that $X \setminus \Sigma$ is regular,
and  $X$ has only {\it terminal singularities} at $\Sigma$, i.e.,
(i) the relative canonical divisor $K_{X/\Delta}$ is $\Bbb Q$-Cartier,
and (ii) for any birational quasi-projective morphism
$\mu : X' \to X$ from a normal
scheme $X'$ and any prime divisor $E$ on $X'$ such that
$\mu (E) \subset \Sigma $,
we have
$$
v_{E}(K_{X'/\Delta} - \mu ^{*}K_{X/\Delta}) > 0
$$
where $v_{E}$ denotes the discrete valuation at the generic point of $E$,

\item"(5)"  any closed fiber $X_s$ is a normal crossing divisor on
$X \setminus \Sigma$,
and the pair $(X, X_s)$ has only {\it log terminal singularities}
at $\Sigma$, i.e.,
for any $\mu : X' \to X$ and $E$ as in (4), we have
$$
v_{E}(K_{X'/\Delta} - \mu ^{*}(K_{X/\Delta} + X_s)) > - 1.
$$
\endroster
{}From \S3 on, we shall assume the additional condition:
\roster
\item"(6)"  the $\omega _{X/\Delta}^{[m]} = \Cal O_X(mK_{X/\Delta})$
satisfy the condition $(S_{3})$ for all $m \in \Bbb Z$.
\endroster
\enddefinition

$f: X \to \Delta$ is said to be {\it minimal} if
$f$ is projective and $K_{X/\Delta}$ is $f$-nef.
The condition (4) actually follows from (5),
since $X_s$ is a Cartier divisor.
By the adjunction of the log canonical divisors,
we deduce from (3) and (5) that the pairs $(S_i, D_i)$ for
$$
D_i = \sum _{j \neq i} S_j \vert _{S_i}
$$
have only log terminal singularities at $\Sigma$.
We note that we excluded the pinch point and its quotients from the
singularities of $X_s$ by the condition (3).
The minimal resolutions of the log terminal singularities of
surfaces are classified by [TM].
In characteristic 0, the condition (6) follows from the condition (4),
since the terminal singularities are Cohen-Macauley in this case.

\definition{1.2. Definition}
In this paper, a {\it curve} is an irreducible and
reduced closed subscheme $C_j$ of a fiber $X_s$
for $s \in \Delta$ which has dimension 1 over the residue field
$\kappa (s) = A_s/\frak m_s$.
As in [KMM], when $f$ is projective, we define
$$\align
N^1(X/\Delta) &= (\{\text{Cartier divisor $L$ on $X$}\}/
\sim _{\text{num}/\Delta})
\otimes \Bbb R \\
N_1(X/\Delta) &= (\{\text{1-cycle $\sum d_jC_j$ in the fibers of
$f$}\}/\sim _{\text{num}/\Delta}) \otimes \Bbb R
\endalign$$
where $\sim _{\text{num}/\Delta}$
denotes the numerical equivalence over $\Delta$.
A $\Bbb Q$-Carier divisor $L$ gives a linear functional $h_L$ on
$N_1(X/\Delta)$ defined by $h_L(z) = (L \cdot z)$.
We also define the {\it Hironaka-Kleiman-Mori cone}
(or the {\it cone of curves})
$\overline {\text{NE}}(X/\Delta)$ to be the closed convex cone
in $N_1(X/\Delta)$ generated by the numerical classes of curves.
\enddefinition

\proclaim{1.3. Cone Theorem}
Under the assumption \rom{(1.1)}, assume moreover that  $f$  is projective.
Then $\overline {\text{NE}}(X/\Delta )$
is locally polyhedral in the half space
defined by $h_{K_{X/\Delta }}(z) < 0$.
Moreover, any extremal ray in this half space
is generated by a class of a curve on a fiber of $f$.
\endproclaim

\demo{Proof}
We have $K_{X/\Delta }\vert _{S_i} = K_{S_i} + D_i$ for any irreducible
component $S_i$ of any fiber of $f$.
By [Mo1] and [TM], the Hironaka-Kleiman-Mori cone
$\overline {\text{NE}}(S_i) \subset  N_1(S_i)$ is locally polyhedral
in the half spaces $h_{K_{S_i} + D_i}(z) < 0$,
and any extremal ray is generated by a curve on $S_i$.

The images of the $\overline {\text{NE}}(S_i)$
by the surjective homomorphism
$$
\bigoplus_i N_1(S_i) \to N_1(X/\Delta )
$$
generate $\overline {\text{NE}}(X/\Delta )$ as a closed convex cone.

Let $X_s$ be a smooth closed fiber and $C$ a curve on
$X_s$ which generates an extremal ray of $\overline {\text{NE}}(X_s)$.
We shall porve that the class of $C$ in $N_1(X/\Delta )$ is contained in
the image of $\overline {\text{NE}}(X_{\eta})$.
Let $\Cal H$  be an irreducible component of the relative Hilbert scheme
of  $f: X \to \Delta$ which contains the point $[C]$ corresponding to $C$.
Let $\bar k$ be an algebraic closure of the residue field  $k = \kappa (s)$,
and $\hat A$ a complete discrete valuation ring
which dominates $A_s$ and whose residue field is $\bar k$.
The base change to $\bar k$ gives a closed subscheme $\bar C$
of $X_{\bar s}$.

By [Mo1, 2.7], (i) $\bar C$ is a disjoint union of $(- 1)$-curves,
(ii) $\bar C$ is a union of fibers of a conic bundle, or
(iii) $\rho (X_s) = 1$.  By an infinitesimal calculus,
the geometric point  $[\bar C] \in \Cal H(\bar k)$ corresponding to
$\bar C$ can be extended to a geometric point
$\hat C \in \Cal H(\hat A)$
in the cases (i) and (ii).
In the case (iii), we have the same assertion if we replace $C$ by a
hyperplane section of $X_s$ belonging to an $f$-very ample linear system.
By the universality of $\Cal H$, we
conclude that there exists a closed subscheme $C_{\eta}$ of $X_{\eta}$
which has the same image as $C$ in $N_1(X/\Delta )$.

Since there are only a finite number of non-smooth fibers of $f$, we have
the theorem. Q.E.D.
\enddemo

We can prove a relative version of the cone theorem,
where $N_1(X/B)$, etc. are defined similarly as in
Definition 1.2 (cf. [KMM]).

\proclaim{1.3'. Theorem}
Let $B$ be a scheme with a flat and quasi-projective morphism
$g : B \to \Delta $, or one which is obtained by a completion from
such a scheme,
and $f' : X \to B$ a projective morphism such that
$f = g \circ f'$ satisfies (1.1).
Then the closed convex cone $\overline {\text{NE}}(X/B)
\subset N_1(X/B)$  is locally polyhedral in the half space
defined by $h_{K_{X/\Delta }}(z) < 0$.
Moreover, any extremal ray in this half space
is generated by a class of a curve on a fiber of $f'$.
\endproclaim

\head 2. Contraction theorem \endhead

\proclaim{2.1. Lemma}
Let $(S, D)$  be a pair of a normal surface
over an algebraically closed field $k$
and a $\Bbb Q$-divisor which has only log terminal
singularities
(thus $\llcorner D \lrcorner = 0$), $L$ a Cartier divisor on $S$,
and let $f: S \to B$ be a projective surjective morphism
whose fibers are at most 1 dimensional.
Assume that
$L - (K_{S} + D)$ is $f$-nef and $f$-big (cf. [KMM]).
Then $R^1f_*\Cal O_S(L) = 0$.
\endproclaim

\demo{Proof}
Let $g: S' \to S$ be the minimal resolution.  We define
$K_{S'} + D' = g^*(K_S + D)$ and $L' = g^*L$.
Then the pair $(S', D')$ and $L'$ satisfy the
conditions of the lemma with respect to $g$ and $f \circ g$.  Therefore,
we may assume that $S$ is nonsingular.
We may also assume that $L - (K_S + D)$ is $f$-ample.

Let  $E = \sum e_iE_i$ be a nonzero effective divisor on $S$
whose irreducible components $E_i$ are contained in a fiber of $f$.
Suppose that $H^1(E, \Cal O_E(L)) \ne 0$. Then by Serre duality,
there is a nonzero element $s \in H^0(E, \Cal O_E(K_S + E - L)).$
If we decrease the $e_i$, we may assume that $s$ generates
$\Cal O_E(K_S + E - L)$ at the generic points of the $E_i$.  Then
$((K_S + E - L)\cdot E_i) \ge 0$, hence
$((E - D).E_i) > 0$, a contradiction.  Q.E.D.
\enddemo

\proclaim{2.2. Lemma}
Let $(S, D)$  be a pair of a normal projective surface
over an algebraically closed field $k$
and a $\Bbb Q$-divisor which has only log terminal
singularities. Assume that
$- (K_{S} + D)$ is ample. Then $H^1(S, {\Cal O}_{S}) = 0$.
\endproclaim

\demo{Proof}
We shall prove that $S$ is a rational surface. Let
$g : S' \to S$ be a minimal resolution.
Since $(S, D)$  has only log terminal singularities,
there exists a $\Bbb Q$-divisor $D'$ on $S'$ such that
\roster
\item"(1)" $K_{S'} + D' - g^*(K_{S} + D)$
is an effective divisor supported on
the exceptional locus of $g$,

\item"(2)" $(S', D')$ has only log terminal singularities,

\item"(3)" $- (K_{S'} + D')$ is ample.
\endroster

Since $\vert mK_{S'}\vert = \emptyset $ for all $m > 0$,
$S'$ is birationally equivalent to a ruled surface.
By (3), $\overline {\text{NE}}(S')$ is spanned by extremal rays,
which are generated by rational curves.  Hence the genus of the base curve
of the ruling is 0.  Q.E.D.
\enddemo

\proclaim{2.3. Contraction Theorem}
Under the assumption of Theorem 1.3, let $R$ be an extremal ray of
$\overline {\text{NE}}(X/\Delta )$ in the half space
$h_{K_{X/\Delta }}(z) < 0$.
Then there exists a projective surjective morphism $\varphi : X \to Z$
to a normal scheme $Z$ which is projective over $\Delta $
such that $\varphi _*\Cal O_X = \Cal O_Z$ and that
any curve  $C$  on any fiber of  $f$  is mapped to a point
by  $\varphi$  if and only if  $\text{cl}(C) \in R$.
\endproclaim

\demo{Proof}
We shall use the contraction theorem of log surfaces (cf. [TM]).
By definition of an extremal ray,
there exists a line bundle  $L$ on $X$ such that
$(L \cdot C) \ge 0$ for any curve $C$ on any fiber of $f$,
and $(L \cdot C) = 0$ if and only if $\text{cl}(C) \in R$.
We shall prove that the linear system $\vert mL\vert$ is relatively free
over $\Delta$ for a positive integer $m$, i.e.,
the canonical homomorphism $f^*f_*\Cal O_X(mL) \to \Cal O_X(mL)$ is
surjective.

Let us fix a closed point $s$ of $\Delta$.
Let $\bar k$ be an algebraic closure of the residue field  $k = \kappa (s)$,
$X_{\bar s}$ the corresponding geometric fiber,
$\bar S_i$ the irreducible components of $X_{\bar s}$, and
$\bar D_i$ the divisors on the $\bar S_i$ obtained from the $D_i$
by the base change.
As in the proof of Theorem 1.3,
let $\hat A$ be a complete discrete valuation ring
which dominates $A_s$ and whose residue field
$\hat A/\frak m_{\bar s} = \bar k$.
Moreover, we assume that $\frak m_{\bar s} = \frak m_s\hat A$.
Let $\hat f: \hat X \to \hat \Delta$ be the morphism obtained by the base
change $A_s \to \hat A$.  By the flat base change theorem,
it is sufficient to prove that
$\hat f^*\hat f_*\Cal O_{\hat X}(m\hat L) \to \Cal O_{\hat X}(m\hat L)$
is surjective.

Let  $\Cal C_{\bar s}$  be set of the curves on $X_{\bar s}$
which are irreducible components of the 1-dimensional
closed subschemes obtained by the base change from the curves $C$
on $X_s$ such that $\text{cl}(C) \in R$.
Then the curves in $\Cal C_{\bar s}$
generates an extremal face $R_{\bar s}$ of
$\overline {\text{NE}}(X_{\bar s})$.
Let $E_{\bar s}$ be the union of these curves.

There are four cases:
\roster
\item"(1)" there are only a finite number of curves in $\Cal C_{\bar s}$,

\item"(2)" $E_{\bar s}$ contains an irreducible component $\bar S_i$
of $X_{\bar s}$ but $E_{\bar s} \ne X_{\bar s}$,

\item"(3)" $E_{\bar s} = X_{\bar s}$ but
$\text{cl}(\bar C) \notin R_{\bar s}$
for some curve $\bar C$ on $X_{\bar s}$,

\item"(4)" $\text{cl}(\bar C) \in R_{\bar s}$ for all
$\bar C$ on $X_{\bar s}$.
\endroster

Case (1).  Let $\bar C$ be a curve $\Cal C_{\bar s}$.
If $\bar C$ is contained in two irreducible components
$\bar S_i$ and $\bar S_j$ of $X_{\bar s}$, then
$(\bar S_i \cdot \bar C) < 0$ and
$(\bar S_j \cdot \bar C) < 0$,
because $\bar C$ is contractible in
the normal surfaces  $\bar S_i$ and $\bar S_j$.  Hence
$\bar C'\subset \bar S_i \cap \bar S_j$ for any
$\bar C' \in \Cal C_{\bar s}$.
Here we note that we used the condition (3) of (1.1).
If $\bar C \subset \bar S_i$,
$\bar C \not\subset \bar D_i$ and $\bar C \cap \bar D_i
\ne \emptyset$,
then $(\bar S_i\cdot \bar C) < 0$,
hence $\bar C'\subset \bar S_i$ for any $\bar C' \in \Cal C_{\bar s}$.
Finally, if $\bar C \subset \bar S_i$ and
$\bar C \cap \bar D_i = \emptyset$, then
$\bar C'\cap \bar D_j = \emptyset$ for any $\bar C' \in \Cal C_{\bar s}$
such that $\bar C'\subset \bar S_j$.
In any of the above cases,
one can prove that there exist a projective scheme $Z_{\bar s}$
of dimension 2
which satisfies the condition  $(S_2)$ and a birational morphism
$\varphi _{\bar s}: X_{\bar s} \to Z_{\bar s}$
whose exceptional locus coincides with  $E_{\bar s}$.

By Lemma 2.1, one can easily check
that  $R^1\varphi_{\bar s*}\Cal O_{X_{\bar s}} = 0$, and
there exists an ample line bundle $H_{\bar s}$ on $Z_{\bar s}$
such that $L_{\bar s} = \varphi _{\bar s}^*H_{\bar s}$ (cf. [A1]).
Therefore, the linear system $\vert mL_{\bar s}\vert$ is free and
$H^1(X_{\bar s}, \Cal O_{X_{\bar s}}(mL_{\bar s})) = 0$
for a positive integer $m$ by the Serre vanishing theorem,
hence our assertion.

If we choose $m$ large enough, then the image $Z$ satisfies conditions
$(S_3)$ and  $(R_2)$, and in particular,
it is normal.

Case (2).  The extremal face  $R_{\bar s}$  determines an extremal face of
$\overline {\text{NE}}(\bar S_i)$ with respect to
$K_{\bar S_i} + \bar D_i$, hence
a contraction morphism
$\varphi _{\bar S_i} : \bar S_i \to \varphi _{\bar S_i}
(\bar S_i) = \bar B$.
We have two subcases:
\roster
\item"(2a)" $\text{dim }\bar B = 1$,

\item"(2b)" $\text{dim }\bar B = 0$.
\endroster

In the subcase (2a), suppose that a general fiber $\bar C$ of
$\varphi _{\bar S_i}$ does
not intersect $\bar D_i$.
Then we have  $(\bar S_j\cdot \bar C')
= 0$  for any $j$ and any $\bar C' \in \Cal C_{\bar s}$.
If $\bar C''$ is an irreducible component of $\bar D_i$,
then $\bar C''$ is contained in a fiber of $\varphi _{\bar S_i}$,
hence $(\bar S_i \cdot \bar C'') = 0$.
So if $\bar C''\subset \bar S_j$, then
the contraction
morphism from  $\bar S_j$ determined by the face $R_{\bar s}$
is of fiber type again, i.e.,
$\bar S_j \subset E_{\bar s}$. Hence
$X_{\bar s} \subset E_{\bar s}$, a contradiction.

Therefore, there is an irreducible component
$\bar S_j$ such that $\bar S_i \cap \bar S_j$
is a section of $\varphi _{\bar S_i}$.
Then  $(\bar S_i \cdot \bar C) < 0$ for a general fiber
$\bar C$ of $\varphi _{\bar S_i}$, and
$(\bar S_i\cdot \bar C') < 0$  for any $\bar C' \in \Cal C_{\bar s}$.
Hence $E_{\bar s} = \bar S_i$.
Let $\bar C''$ be an irreducible component of
$\bar D_i$ in a fiber of
$\varphi _{\bar S_i}$.
Since $(\bar S_j \cdot \bar C'') > 0$,
$\bar C''$ intersects $\bar S_j$
at a nonsingular point of $\bar S_i$.
In particular, $\bar C''$
is the only component of  $\bar D_i$ in the same fiber as
$\bar C''$.

Then one can construct a projective scheme $Z_{\bar s}$ of dimension 2
which satisfies the condition  $(S_2)$ and a morphism
$\varphi _{\bar s}: X_{\bar s} \to Z_{\bar s}$
which coincides with  $\varphi _{\bar S_i}$ on $\bar S_i$
and is isomorphic outside  $\bar S_i$.
As in the case (1), we have only to prove that
$R^1\varphi _{\bar s*}\Cal O_{X_{\bar s}} = 0$.

Let  $X'_{\bar s}$  be the union of the irreducible components
of  $X_{\bar s}$ other than  $\bar S_i$.
Thus $\bar D_i = X'_{\bar s} \cap \bar S_i$.
We have exact sequences
$$\align
0 &\to \Cal O_{X'_{\bar s}}(- \bar D_i) \to
\Cal O_{X_{\bar s}} \to \Cal O_{\bar S_i} \to 0 \\
0 &\to \Cal O_{X'_{\bar s}}(- \bar D_i) \to
\Cal O_{X'_{\bar s}} \to \Cal O_{\bar D_i} \to 0.
\endalign$$
We have $R^1\varphi _{\bar S_i*} \Cal O_{\bar S_i} = 0$ by Lemma 2.1.
Since $R^1\varphi _{\bar s*} \Cal O_{X'_{\bar s}} = 0$.
and $\bar D_i$  is connected, we have also
$R^1\varphi _{\bar s*} \Cal O_{X'_{\bar s}}(- \bar D_i) = 0$.
Hence $R^1\varphi _{\bar s*}\Cal O_{X_{\bar s}} = 0$.

In the subcase (2b), since  $(\bar S_i\cdot \bar C) < 0$
for any $\bar C \in \Cal C_{\bar s}$,
we have again $E_{\bar s} = \bar S_i$.
Moreover, since any irreducible component of $\bar D_i$
is an ample $\Bbb Q$-Cartier divisor
on $\bar S_i$, it intersects any curve on $\bar S_i$.
But since $K_{\bar S_i} + \bar D_i$
is negative, $\bar D_i$ consists of 1 or 2 irreducible
components.  Then Lemma 2.2 combined with the argument
in the subcase (2a) yields the result.

Case (3).  Let  $\bar S_i$ be an irreducible component of $X_{\bar s}$
such that
the image of the associated contraction morphism $\varphi _{\bar S_i}$
is a curve.  If  $\bar D_i$ contains a section of
$\varphi _{\bar S_i}$, then
we have $(\bar S_i \cdot \bar C) < 0$ for $\bar C \in \Cal C_{\bar s}$,
and $E_{\bar s} = \bar S_i$ as in the case (2), a contradiction.
Hence $(\bar S_i \cdot \bar C) = 0$ for any $\bar C \in \Cal C_{\bar s}$,
and no irreducible component of
$X_{\bar s}$ is contracted to a point by the associated contraction morphism.
Any irreducible component of any intersection
$\bar S_i \cap \bar S_j$ is at the same time a
fiber of  $\varphi _{\bar S_i}$  and $\varphi _{\bar S_j}$.
Thus we have a
1-dimensional reduced scheme $Z_{\bar s}$ with only ordinary double points
as singularities
and a contraction morphism  $\varphi _{\bar s}:X_{\bar s} \to Z_{\bar s}$
with $R^1\varphi _{\bar s*}\Cal O_{X_{\bar s}} = 0$, and we are done.

Case (4).  We have just  $\varphi = f$.
Q.E.D.
\enddemo

We have again a relative version of the contraction theorem.

\proclaim{2.3'. Theorem}
Under the assumption of Theorem 1.3', if  $R$  is an extremal ray of
$\overline {\text{NE}}(X/B)$
in the half space $h_{K_{X/\Delta }}(z) < 0$,
then there exists a projective surjective morphism
$\varphi : X \to Z$ to a normal scheme  $Z$ which is projective over $B$
such that $\varphi _*\Cal O_X = \Cal O_Z$ and that
any curve  $C$  on any fiber of  $f'$  is mapped to a point
by  $\varphi$  if and only if  $\text{cl}(C) \in R$.
\endproclaim

\head 3. Canonical cover \endhead

We shall extend the construction of the canonical cover (cf. [R]).
We assume that the morphism $f: X \to \Delta$ satisfies (1.1) including
the condition (6).
Moreover, we assume that $A$ is a complete discrete valuation ring whose
residue field $k$ is algebraically closed.
For a closed point $P \in \Sigma$, the positive integer
$$
r = \min \{m \in \Bbb N \: \text{$mK_{X/\Delta}$ is Cartier at $P$}\}
$$
is called the {\it index} of $K_{X/\Delta}$ at $P$.
We note that $K_{X/\Delta}$ is Cartier on $X \setminus \Sigma$.

Let us fix a point $P \in \Sigma$.
By replacing $X$ by a small neighborhood of $P$ if necessary,
we assume that $K_{X/\Delta}$ is Cartier on $U = X\setminus \{P\}$
and there exists
a nowhere vanishing section $\theta $ of $\Cal O_X(rK_{X/\Delta})$.
Then an $\Cal O_X$-algebra structure on
$$
\Cal F = \bigoplus _{m=0}^{r-1}\Cal O_X(- mK_{X/\Delta})t^m
$$
is defined by an equation $t^r = \theta $.
Let $Y = \text{Spec }\Cal F$, $\pi : Y \to X$ the natural projection,
and $Q = \pi ^{- 1}(P)$.
By the condition (6) of (1.1), $Y$ satisfies the condition $(S_{3})$.

The group scheme $\mu _{r} = \text{Spec }\Bbb Z[\zeta ]/(\zeta ^r - 1)$
with co-multiplication
$\zeta \mapsto \zeta \otimes \zeta $ acts on $Y$ by the co-action
$$
x \mapsto \zeta ^m \otimes x \text{ if $x \in \Cal O_X(- mK_{X/\Delta})t^m$}.
$$
The ring of $\mu _r$-invariants of $\Cal O_Y$ coincides with $\Cal O_X$.

Since the choice of $\theta $ is large, the {\it canonical cover}
$\pi : Y \to X$
is not unique even locally in the etale topology.

\proclaim{3.1. Theorem} If $X$ is a small enough neighborhood of $P$
and $\theta $ is chosen sufficiently general, then
$f\circ \pi : Y \to \Delta$ satisfies the conditions of (1.1).
Moreover, $Y$ is Gorenstein.
\endproclaim

\demo{Proof}
By [H],
$$
\pi _*\omega _{Y/\Delta} \simeq \Cal Hom_{\Cal O_X}(\pi _*\Cal O_Y,
\omega _{X/\Delta})
= \bigoplus _{m=0}^{r-1}\Cal O_X((m + 1)K_{X/\Delta})t^{-m}
\simeq  \pi _*\Cal O_Y \cdot t
$$
as $\pi _*\Cal O_Y$-modules.
Hence $Y$ is Gorenstein.

One can decompose $\pi$ as follows.
Let  $p$  be the characteristic of $k$,
and let $r = \tilde rp^e$ with $(\tilde r, p) = 1$.
We set $\tilde t = t^{p^e}$, and
$\tilde Y = \text{Spec }(\bigoplus _{m=0}^{\tilde r-1}
\Cal O_X(- mp^eK_{X/\Delta})\tilde t^m)$ with $\tilde t^{\tilde r}
= \theta$.
Let $\tilde \pi : \tilde Y \to X$ and $\sigma : Y \to \tilde Y$
be the natural morphisms.
Then $\tilde \pi $ is etale over $X \setminus \Sigma$ and $\sigma$ is purely
inseparable.

Let $T_i$ be the inverse image by $\pi$ of an irreducible component $S_i$
of $X_s$.  If $T_i$ is not irreducible, then its irreducible components
meet only at $Q$. But there is a positive integer $n$
such that $nT_i$ is a Cartier divisor on $Y$ and satisfies the
condition $(S_2)$, a contradiction.
Hence $T_i$ is an irreducible $\Bbb Q$-Cartier divisor.

We shall prove that $Y$ is normal.
If $r = \tilde r$, then it is clear that $Y \setminus \{Q\}$ is regular,
hence $Y$ is normal.
We assume that $p\vert r$ in the following.
For a point $P' \ne P$ of $X$,
let $\theta _0$ be a germ of a generating section of
$\Cal O_X(K_{X/\Delta})$ at $P'$,
and write $\theta = h\theta _0^r$ for a germ of an invertible function
$h$ at $P'$.
Let $s$ be a nonzero global section of $\Cal O_X(K_{X/\Delta})$.
If we write $s = c\theta _0$ for a germ of a regular function $c$ at $P'$,
then $s^r = c^r\theta _0^r = g\theta$ for a global regular function $g$.
Thus a local section $dh$ of $\Omega ^1_{X/\Delta}$ vanishes if and only if
$dg/g$ vanishes.
If $U'$ is the open subset of $X$
on which $f$ is smooth and $dg/g$ does not vanish, then
$f\circ \pi$ is smooth over $\pi ^{- 1}(U')$.
If we replace  $\theta$ by $u\theta$ for an invertible function $u$, then
$dg/g$ is replaced by $dg/g - du/u$.
By Bertini's theorem, $du/u$ has no zero on $U$ for general
$u$.  The 1-form $dg/g$ gives a rational section of the projection
$\Bbb P_*(\Omega ^1_{U/\Delta}) \to U$,
and its intersection
with the section given by $du/u$ is at most 2-dimensional for general $u$,
since the latter is free.  Then one can arrange $u$ so that
the zero locus of $dg/g - du/u$ is at most 1-dimensional.
Therefore, $X \setminus U'$ is at most 1-dimensional, and $Y$ is normal.
The above argument also showed that the condition (2) of (1.1) holds for $Y$.

We shall prove that $Y$ has only terminal singularities.
Then it follows the conditions (1) and (4) of (1.1) for $Y$.
Let $\mu :X'\to X$ be a quasi-projective birational morphism from a normal
scheme $X'$, and $E$ a prime divisor on $X'$
such that $\text{codim }\mu (E) \ge 2$.
We have $K_{X'/\Delta} = \mu ^*K_{X/\Delta} + \Gamma $
for some effective $\Bbb Q$-divisor $\Gamma $ whose
support contains $E$; we have $\Gamma = \frac{a}{r}E + ..$.
with $a \in \Bbb N$. We define
an $\Cal O_{X'}$-algebra structure on
$$
\Cal F' = \bigoplus _{m=0}^{r-1}\Cal O_{X'}
(- mK_{X'/\Delta} + \llcorner m\Gamma \lrcorner )t^m
$$
by $t^r = \mu ^*\theta \in \Cal O_{X'}(rK_{X'/\Delta} - r\Gamma )$.
Let $Y' = \text{Spec }\Cal F'$, and $\pi ': Y' \to X'$
the projection.
Since there are natural homomorphisms
$$
\mu ^*\Cal O_X(-mK_{X/\Delta}) \to \Cal O_{X'}
(- mK_{X'/\Delta} + \llcorner m\Gamma \lrcorner )
$$
we have a birational morphism $\mu ':Y'\to Y$ which covers $\mu $.

We have to prove that $Y'$ is normal.
By construction, $Y'$ satisfies
the condition $(S_2)$.
Since $\mu '$ is an isomorphism outside the support
of $\Gamma $, it is enough to prove that $Y'$ is regular above the generic
point $\eta _E$ of $E$.
Let $v$ and $\theta '_0$ be generators of the maximal ideal
$\frak m_{E,X'}$ and $\omega _{X'/\Delta}$ at $\eta _E$, respectively.
We can write $\mu ^*\theta = hv^a\theta _0^{\prime r}$ for
an invertible function $h$ at $\eta _E$.
Let $c = (r, a)$, $r = cr'$ and $a = ca'$.
Then $am_0 - (\ulcorner am_0/r' \urcorner - 1)r' = 1$
for an integer $m_0$ such that $0 < m_0 \le r'$.
We set
$$\align
u_j &= \theta _0^{\prime -jr'}v^{-ja'}t^{jr'}
\in \Cal O_{X',\eta _E}(- jr'K_{X'/\Delta} +
\llcorner jr'\Gamma \lrcorner )t^{jr'}
\text{ for $0 \le j < c$} \\
v' &= \theta _0^{\prime -cm_0}v^{- \ulcorner am_0/r' \urcorner + 1}t^{cm_0}
\in \Cal O_{X',\eta _E}(- cm_0K_{X'/\Delta} +
\llcorner cm_0\Gamma \lrcorner)t^{cm_0}.
\endalign
$$
Then $u_j^c = h^j \in \Cal O_{X',\eta _E}$ is invertible,
and $v'{}^{r'} = h^{m_0}v$.
Hence $v'$ generates the maximal ideals at the points of $Y'$
above $\eta _E$, and $Y'$ is regular there.
Therefore, $Y'$ is isomorphic to the normalization of the
fiber product $Y \times _X X'$.

We have
$$
\pi '_*\omega _{Y'/\Delta} \simeq \bigoplus _{m=0}^{r-1}\Cal O_{X'}
((m+1)K_{X'/\Delta}- \llcorner m\Gamma \lrcorner )t^{-m},
$$
and there are natural homomorphisms
$$
\mu ^*\Cal O_X((m + 1)K_{X/\Delta}) \to
\Cal O_{X'}((m + 1)K_{X'/\Delta} - \llcorner m\Gamma \lrcorner )
$$
which is not isomorphic at $\eta _E$ for $m = r - 1$.
Hence the discrepancy
$K_{Y'/\Delta} - \mu '{}^*K_{Y/\Delta}$
has positive coefficients at prime divisors of $Y'$
which lie above $E$.
Thus $Y$ has only terminal singularities.

Next, we shall prove that the pair $(Y, Y_s)$ has only a log terminal
singularity at $Q$. Supposing that $\mu (E) = P$,
let us write $\mu ^*X_s = bE + ... $ for $b \in \Bbb N$.  Since $(X, X_s)$
is log terminal at $P$, we have $a/r - b > - 1$.  The coefficient of the
discrepancy $K_{Y'/\Delta} - \mu '{}^*(K_{Y/\Delta} + Y_s)$ at the points
above $\eta _E$ is
$r'(a - \llcorner a(r - 1)/r \lrcorner - b)
= r'(\ulcorner a/r \urcorner - b) \ge 0$.
Hence $(Y, Y_s)$ is log terminal at $Q$.

Let $T'_i$ be the normalizations of the $T_i$, and
$F_i = \sum _{j \neq i} T_j \vert _{T'_i}$.
Since the pair $(Y, Y_s)$ is log terminal at $Q$,
so are the pairs $(T'_i, F_i)$ above $Q$ by the adjunction.
Then it follows easily that the $F_i$ are reduced, hence
$Y_s$ is a normal crossing divisor on $Y \setminus \{Q\}$.
The restriction of a local section of $\Cal O_{T_i}$
on $T_i \setminus \{Q\}$ to $F_i \setminus \{Q\}$ extends
to sections of the $\Cal O_{T_j}$ on the $T_j \setminus \{Q\}$ for
$j \ne i$, hence a section of $\Cal O_{Y_s}$ over $Y_s \setminus \{Q\}$.
Since $Y_s$ satisfies the condition $(S_2)$,
it extends over to $Y_s$, hence the $T_i$ are normal.
We note that this is the point where we needed the condition (6) of (1.1).
This completes the proof of the theorem.
Q.E.D.
\enddemo

\head 4. Classification of singularities \endhead

Let $f: X \to \Delta $ be as in \S3. We shall classify the
singular points of $X$ (cf. [Mo2]).

\proclaim{4.1. Theorem}
Let $P \in X$ be a point of index $r > 1$, and $\pi : Y \to X$, etc.
as in \S3. Then the completion of $\Cal O_{Y,Q}$ at $Q = \pi ^{-1}(P)$
with the $\mu _r$-action is isomorphic to the completion of
$A[x, y, z]/(F)$
for some semi-invariant coordinates $x$, $y$ and $z$, where
$$
F = xy + G(z^r)
$$
for some polynomial $G \in A[z^r]$,
and the action of $\mu _r$ is given by
$$
x \mapsto \zeta ^a \otimes x, \quad y \mapsto \zeta ^{-a} \otimes y,
\quad z \mapsto \zeta \otimes z
$$
for some positive integer $a$ such that $(r, a) = 1$. Moreover, if
$\Cal O_{Y_s,Q}$ is not integral, then
$$
F = xy + \tau
$$
for a generator $\tau$ of the maximal ideal of $A$.
\endproclaim

\demo{Proof}
By the classification of log terminal singularities
of surfaces of index 1,
the completion of $\Cal O_{Y_s,Q}$ with the $\mu _r$-action
is isomorphic
to a hypersurface singularity $k[[x, y, z]]/(F_s)$ for some
semi-invariant coordinates $x$, $y$ and $z$, where one of the following
holds (cf. [A2]):
\roster
\item"(1)" $F_s = xy$,

\item"(2)" $F_s$ has at most an isolated singularity
and $\text{ord}(F_s) \le 2$.
\endroster

Since $Y_s = X_s \times _X Y$,
the completion of $\Cal O_{Y,Q}$ with the $\mu _r$-action
is isomorphic to
the completion of $A[x, y, z]/(F)$ such that $F \equiv F_s
(\text{mod } \tau A[x, y, z])$.
We note that the elements in $A$ are $\mu _r$-invariants.
Since the action
of $\mu _r$ is nontrivial and $Y$ has only isolated singularity, the
subspace of $\text{Spec }A[x, y, z]$ defined by $F = 0$
intersects the fixed
locus of $\mu _r$ only at the origin.
In fact, the fixed locus of codimension two or
more has toric singularities, while if there is only one
coordinate which is not invariant, then $\omega _{X/\Delta,P}$
is generated by one element over $\Cal O_{X,P}$, a contradiction.
Therefore, none of $x$, $y$ nor $z$
is $\mu _{r'}$-invariant for an integer $r'$ with $r'\vert r$ and $r'> 1$.
Moreover, $F$ contains a constant term on $x, y$ and $z$,
hence is a $\mu _r$-invariant.
Then so is $F_s$.

Thus we have the first assertion of the theorem in the case (1).
We write $F = xy + \tau G'(z^r)$.
Since $Y$ has only isolated singularities, $G'$ is not divisible by $\tau$.
Since the $S_i$ are $\Bbb Q$-Cartier, the prime divisor
on $Y$ defined by $x = \tau = 0$ should be $\Bbb Q$-Cartier. Then $G'$ is
invertible, hence the last assertion.

In the case (2), if $\text{ord}(F_s) = 1$, then there is an invariant
coordinate, a contradiction. Hence $\text{ord}(F_s) = 2$.
If $F_s$ contains a
degree 2 term of the form $xy$, then we are done. If it contains
$x^2$ and there are no other terms of degree 2, then $r = 2$ and there
is a degree 3 term, a contradiction to the fact that none of the coordinates
$x$, $y$ and $z$ are $\mu _r$-invariants. Q.E.D.
\enddemo

\proclaim{4.2. Corollary}
There exists a projective birational morphism
$\mu :X'\to X$ from a normal scheme $X'$
such that the discrepancy coefficient of one of the exceptional divisors
with respect to $K_{X/\Delta}$ is equal to $1/r$.
\endproclaim

\demo{Proof}
We take a weighted blowing up with weights $\frac{1}{r}(a, r - a, 1, r)$
for the semi-invariants $(x, y, z, \tau )$ (cf. [KKMS] and [K4]). Q.E.D.
\enddemo

In the above theorem, the {\it axial multiplicity} of the singular point
$P \in X$ is defined to be largest integer $n$
such that $G (\text{mod }z^rA[z^r]) \in \tau ^nA$.
If $G_s \equiv G (\text{mod }\tau A[z^r]) \in z^rk[z^r]$
is equal to 0 or contains a linear term in $z^r$, then the
singular point $P \in X$ is said to be of {\it simple type}.
Then we call the above $\mu :X' \to X$ the {\it standard blowing up} of $X$.
In this case, the exceptional locus $E$ of $\mu$ is an irreducible divisor
with $\rho (E) = 1$, and
$X'$ satisfies the conditions in (1.1) and
has at most 3 terminal singular points which are of simple type again;
it has singular points of indices $a$ and $r - a$ if $a > 1$ and
$r - a > 1$, respectively, and
a singular point of index $r$ and axial multiplicity $n - 1$ if $n > 1$.

\proclaim{4.3. Corollary}
Let $S_i$ be an irreducible component of the closed fiber $X_s$ of $f$, and
let $D_i$ be as before.  Then the singularities of the pair $(S_i, D_i)$
are toric;
a point on $S_i \setminus D_i$ is of type
$\frac{1}{hr^2}(bhr-1, 1)$,
while $S_i$ is of type $\frac{1}{r}(b, 1)$ at a point on $D_i$,
for some positive integers $r$, $b$ and $h$ such that
$0 < b < r$ and $(r,b)=1$.
Moreover, the former is of simple type if and only if  $h = 1$,
and the latter is always of simple type.
\endproclaim

\demo{Proof}
In the former case, we have  $F_s = xy + z^{hr}$
for some positive integer  $h$.
Then we set  $ab = cr + 1$ for a positive integer $c$.  Q.E.D.
\enddemo

\proclaim{4.4. Theorem}
Let $P \in X$ be a point of index 1. If $f$ is not
smooth at $P$, then the completion of $\Cal O_{X,P}$ is isomorphic to the
completion of $A[x, y, z]/(F)$ for some coordinates $x$, $y$ and $z$,
where one of the following holds:
\roster
\item"(1)" $F = xyz + \tau $,

\item"(2)" $F = xy + \tau $,

\item"(3)" $F \equiv F_s (\text{mod }\tau A[x, y, z])$ for a polynomial
$F_s \in k[x, y, z]$
which defines a rational double point.
\endroster
\endproclaim

\demo{Proof}
By the classification of log terminal singularities of surfaces of
index 1, we have either $F_s = xyz$, $F_s = xy$, or the case (3).
Since the $S_i$
are $\Bbb Q$-Cartier, we have (1) or (2) in the former cases
as in the proof of Theorem 4.1.  Q.E.D.
\enddemo

\head 5. Flip theorem \endhead

\definition{5.1. Definition}
A {\it flip} in this paper is a diagram
$$
X @>\varphi >>Z @<\varphi ^+<< X^+
$$
which satisfies the following conditions:
\roster
\item"(1)" $\varphi :X \to Z$ and $\varphi ^+: X^+ \to Z$ are
projective birational morphisms,

\item"(2)" $Z$ is a normal scheme with a flat quasi-projective morphism
$g: Z \to \Delta$,
or one which is obtained by a completion from such a scheme,
and the induced
morphisms $f: X \to \Delta $ and $f^+: X^+ \to \Delta$ satisfy
the conditions in \rom{(1.1)},

\item"(3)" $\varphi _{\eta }: X_{\eta } \to Z_{\eta }$ and
$\varphi ^+_{\eta }: X^+_{\eta } \to Z_{\eta }$  are isomorphisms,

\item"(4)" $\varphi $ and $\varphi ^+$ contract only finitely many
curves on finite numbers of closed fibers
of $f$ and $f^+$, respectively,

\item"(5)" $- K_{X/\Delta}$ is $\varphi $-ample and $K_{X^+/\Delta}$
is $\varphi ^+$-ample.
\endroster
\enddefinition

\proclaim{5.2. Theorem}
Let $X @>\varphi >>Z @<\varphi ^+<< X^+$ be a flip, and let $r$ (resp.
$r^+$) be the maximun of the indices of $K_{X/\Delta}$
(resp. $K_{X^+/\Delta}$) at the points
on the exceptional locus $\text{Exc}(\varphi )$
(resp. $\text{Exc}(\varphi ^+)$). Then $r > r^+ \ge 1$.
In particular, the termination of flips holds, i.e.,
there exists no infinite sequence of flips over $\Delta $ as
follows:
$$
X^{(0)} @>\varphi ^{(0)}>> Z^{(0)} @<\varphi ^{(0)+}<< X^{(1)}
@>\varphi ^{(1)}>> Z^{(1)} @<\varphi ^{(1)+}<< X^{(2)} @>\varphi ^{(2)}>>...
$$
\endproclaim

\demo{Proof}
We may assume that $A$ is a complete discrete valuation ring whose residue
field is algebraically closed.
First, we shall prove that $r > 1$. Suppose that $r = 1$.
Let $C$ be an irreducible component of $\text{Exc}(\varphi )$.
Then there is only one irreducible component $S_i$
of the closed fiber $X_s$ which contains $C$.
In fact, if  $C$ is contained in another
component $S_j$, then we have $(S_i \cdot C) \le - 1$ and
$(S_j \cdot C) \le - 1$ by Theorem 4.4.  Since $(X_s\cdot C) = 0$,
there are at least 2 other components of $X_s$ which intersect
$C$ transversally.  But then we have
$(K_{X/\Delta}\cdot C) = ((K_{S_i} + D_i)\cdot C) \ge 0$,
a contradiction. Thus $C \not\subset D_i$.

Since the strict
transform of $C$ in the minimal resolution of $S_i$ is a $(- 1)$-curve,
if $C \cap D_i \ne \emptyset$, then
$(K_{X/\Delta}\cdot C) \ge 0$.
Thus $C \cap D_i = \emptyset$.
The contractibility
of $C$ in $S_i$ implies that there is at most one singular point of $S_i$
on $C$, and it is of type A. Then the surface $\varphi (S_i)$ is smooth at
$Q = \varphi (C)$,
and so is $Z$ at $Q$. But $K_Z$ must not be a $\Bbb Q$-Cartier divisor at
$Q$, a contradiction, hence $r > 1$.

If $r^+ > 1$, then there exists an exceptional divisor $E$ over $X^+$
whose discrepancy coefficient with respect to
$K_{X^+/\Delta}$ is equal to $1/r^+$.
Then the discrepancy coefficient of $E$ with respect to
$K_{X/\Delta}$ is less than $1/r^+$, hence $1/r < 1/r^+$. Q.E.D.
\enddemo

\proclaim{5.3. Theorem}
If the exceptional locus $E$ of the contraction morphism
$\varphi : X \to Z$  in Theorems 2.3 or 2.3'
coincides with an irreducible component  $S_i$ of a closed fiber $X_s$,
then the induced morphism $g: Z \to \Delta $ satisfies the conditions of
\rom{(1.1)}.  If  $X$  is  $\Bbb Q$-factorial, and if $E$ contains a prime
divisor of  $X$ but is not equal to the whole $X$,
then  $E$ coincides with the prime divisor,
$g: Z \to \Delta $ satisfies the conditions of \rom{(1.1)},
and  $Z$  is again $\Bbb Q$-factorial.
\endproclaim

\demo{Proof}
All statements are clear by construction
except the condition (6) of (1.1) (cf. [KMM]).
We may assume that $A$ is a complete discrete valuation ring
whose residue field $k$ is algebraically closed, because the depth is
preserved by the base change to such $A$, but we lose the condition
that $\rho (X/Z) = 1$.
We fix a positive integer $m$.

First, we consider the case in which the closed fiber $E_s$ of $E$ is
1-dimensional in the latter half of the theorem.
We write $K_{X/\Delta} = \varphi ^*K_{Z/\Delta} + dE$
for a positive rational number $d$.
Then we have
$R^1\varphi _*\Cal O_X(mK_{X/\Delta} - \ulcorner mdE \urcorner) = 0$
by Lemma 2.1 as in the proof of Theorem 2.3.
Since $\varphi _{s*}\Cal O_{X_s}(mK_{X_s} - \ulcorner mdE_s \urcorner)
= \Cal O_{Z_s}(mK_{Z_s})$
satisfies the condition $(S_2)$,
$\varphi _*\Cal O_X(mK_{X/\Delta} - \ulcorner mdE \urcorner)
= \Cal O_Z(mK_{Z/\Delta})$
satisfies $(S_3)$.

In the rest of the proof, we assume that $E$ is an irreducible
component of $X_s$.
We write $S = S_i = E$ and $D = D_i$.
For any Weil divisor $W$ on $X$, we shall write
$\Cal O_D(W) = (\Cal O_X(W) \otimes \Cal O_D)/\text{torsion}$
in this proof.
Let $d$ be the smallest positive integer such that
$\Cal O_D(mK_{X/\Delta} - dS)$
has a nonzero global section.
We set $L = mK_{X/\Delta} - dS$.
Let $X'_s$ be the union of the irreducible components of $X_s$
except $S$.
We have an exact sequence
$$
0 \to \Cal O_S(L + S) \to \Cal O_{X_s}(L)
\to \Cal O_{X'_s}(L) \to 0.
$$
We shall prove
\roster
\item"(1)" $R^1\varphi _*\Cal O_S(L + S) = 0$.

\item"(2)" $R^1\varphi _*\Cal O_{X'_s}(L) = 0$.
\endroster

Then (1) implies that the natural homomorphism
$\varphi _*\Cal O_{X_s}(L) \to \varphi _*\Cal O_{X'_s}(L)$
is surjective, hence
$\varphi _*\Cal O_{X_s}(L)$ satisfies the condition
$(S_2)$.  By the results in \S4, $\Cal O_X(L)$ satisfies $(S_3)$,
and the natural homomorphism
$\Cal O_X(L) \to \Cal O_{X_s}(L)$ is surjective.
Then by (2), we have
$R^1\varphi _*\Cal O_{X_s}(L) = 0$,
and the natural homomorphism
$\Cal O_Z(mK_{Z/\Delta}) = \varphi _*\Cal O_X(L) \to
\varphi _*\Cal O_{X_s}(L)$ is surjective, hence
$\Cal O_Z(mK_{Z/\Delta})$ satisfies $(S_3)$.

Our assertions follow immediately from Lemma 2.1
except (1) when $\varphi (S)$ is a point.
There are 3 cases:
\roster
\item"(i)" $D$ consists of 2 rational curves,

\item"(ii)" $D$ is an irreducible rational curve and there are at most
2 singular points of $S$ on $D$,

\item"(iii)" $D$ is an irreducible rational curve and there are 3
singular points of $S$ on $D$.
\endroster

We consider the cases (i) and (ii).  Then
$\Cal O_D(- K_{X/\Delta})$ has a nonzero global section,
hence so has the $\Cal O_D(L - nK_{X/\Delta})$ for any positive integer $n$.
If $n$ is sufficiently large, then we have
$H^1(S, \Cal O_S(L - D - nK_{X/\Delta})) = 0$ by the Serre vanishing theorem.
By the exact sequence
$$
0 \to \Cal O_S(L - D - nK_{X/\Delta}) \to \Cal O_S(L - nK_{X/\Delta})
\to \Cal O_D(L - nK_{X/\Delta}) \to 0,
$$
we have
$H^0(S, \Cal O_S(L - nK_{X/\Delta})) \ne 0$
if $H^1(S, \Cal O_S(L - D - nK_{X/\Delta})) = 0$.
Let $M_n$ be an effective Weil divisor on $S$ which corresponds to a
nonzero global section of $\Cal O_S(L - nK_{X/\Delta})$.
Then it is an ample $\Bbb Q$-Cartier divisor.
By [LM, Lemmas 1 and 2], where the same arguments can be extended
to the $\Bbb Q$-divisors, we have $H^0(M_n, \Cal O_{M_n}) = k$.
Since $H^1(S, \Cal O_S) = 0$ by Lemma 2.2, we have
$H^1(S, \Cal O_S(- M_n)) \allowmathbreak= 0$, hence
$H^1(S, \Cal O_S(L - D -  (n - 1)K_{X/\Delta})) = 0$.
By the descending induction on $n$, we obtain
$H^1(S, \Cal O_S(L - D)) = 0$, i.e., the assertion (1).

In the case (iii), if $\rho (S) = 1$, then we have our theorem
by the following lemma.
We shall finish the proof in the case $\rho (S) > 1$ by using Theorem 5.5.
Here we note that the proof of Theorem 5.5 uses Theorem 5.3 in the case
in which $\rho (S_i/\varphi (S_i)) = 1$.

We may replace $Z$ by its completion at $\varphi (S)$.
Then $\rho (X/Z) = \rho (S)$. In fact, we have
$R^2\varphi _{s*}\Cal O_{X_s} = 0$ by Lemma 2.2, hence
$R^2\varphi _*\Cal O_X = 0$.
By an exact sequence
$$
0 \to \Cal O_X(- X_s) \to \Cal O^*_X \to \Cal O^*_{X_s} \to 0
$$
the natural homomorphism $H^1(X, \Cal O^*_X) \to H^1(X_s, \Cal O^*_{X_s})$
is surjective. Since \linebreak
$H^1(X_s, \Cal O^*_{X_s}) \to H^1(S, \Cal O^*_S)$ is also surjective,
we have our claim.

We apply the minimal model program to the morphism $\varphi :X \to Z$.
By Theorems 5.5 and 5.2, there exists a sequence of flips
$X = X_0 \dasharrow X_1 \dasharrow ...\dasharrow X_t$
and a divisorial contraction $X_t \to Z'$.
By the previous arguments, $Z' \to \Delta$ satisfies the conditions (1.1).
Since the induced morphism
$\psi : Z' \to Z$ contracts only a finite number of curves, we have
$K_{Z'/\Delta} = \psi ^*K_{Z/\Delta}$.
By Lemma 2.1, we have $R^1\psi _*\Cal O_{Z'}(mK_{Z'/\Delta}) = 0$
for any integer $m$, hence
$\psi _*\Cal O_{Z'}(mK_{Z'/\Delta}) = \Cal O_Z(mK_{Z/\Delta})$ satisfies
(1.1).  Q.E.D.
\enddemo

\proclaim{5.4. Lemma}
In the above theorem, assume that $A$ is a complete discrete valuation ring
with an algebraically closed residue field,
$\varphi (S_i)$ is a point,
$D_i$ is irreducible,
$\rho (S_i) = \text{dim } N_1(S_i) = 1$, and that
there are 3 singular points of $S_i$ on $D_i$.
Then the conclusion of the theorem holds.
\endproclaim

\demo{Proof}
We shall use the arguments of [MT].
We explain which part we should pick up from [MT].
We note that the characteristic is assumed to be zero in [MT], but
the arguments which we pick up are characteristic free.

We use the following notation of [MT];
$\bar V = S_i$, $\bar E = D_i$, $f: V \to \bar V$
is the minimal resolution, $E$ is the strict transform of $\bar E$ by
$f$, $D$ is the sum of $E$ and the exceptional divisors of $f$,
$D = \sum_i D_i$ is the decomposition into irreducible components,
$K_V + D^{\#} = f^*(K_{\bar V} + \bar E)$,
$T_i$ for $i = 1, 2, 3$ are exceptional divisors of $f$ above the
3 singular points of $\bar V$ on $\bar E$,
$C_i$ is the irreducible component of $T_i$ which meets $E$,
$C_1 = T_1$, $(C_1^2) = - 2$,
and $F = E + T_1 + T_2 + T_3$.

A {\it permissible linear pencil of rational curves}
$\Lambda$ on $V$ is defined on [MT, p.274].
$L_i$ for $i = 1, 2, 3$ are singular members of $\Lambda$,
and $L_1 = C_1 + 2\ell + D_1$.
$E$ becomes a cross section of $\Lambda$ ([MT, Lemma 2.2]).
Since $(f_*(\ell)^2) > 0$, $D_1$ intersects some irreducible components
$D_2, ..., D_{e+1}$ of $D - D_1$.
They are also cross sections of $\Lambda$.
By [MT, 2.3.1], The $L_i$ are the only singular fibers of $\Lambda$, and
they consist of $(- 1)$-curves and irreducible components of $D$.
If $e = 1$, then there are exactly one $(- 1)$-curve on each $L_i$
by [MT, 2.3.2], and they do not intersect the section $D_2$.
Thus the connected
component of $D$ containing $D_2$ cannot be contracted to a
singular point described in Corollary 4.3,
because its dual graph is a fork, a contradiction.
If $e > 2$, then we have a fork with the central component $D_1$.
If $e = 2$, then there exists a singular fiber $L_2$ which has only
one  $(- 1)$-curve by [MT, 2.3.2].  Thus $D_2$ and $D_3$ should intersect
an irreducible component of $L_2$ which is contained in $D$.
Hence we have a cycle in $D$, a contradiction.
Therefore, it is sufficient to prove the existence of a
permissible linear pencil of rational curves in the following.

In [MT, \S3], by contracting the curves in $D$ except $E$ and $C_1$ on $V$,
we obtain a morphism $g$ to factorize $f$ as $V @>g>> W @>g'>> \bar V$.
Let $\tilde C_1 = g_*(C_1)$ and $\tilde E = g_*(E)$.
Since $\rho (W) = 2$ and $(\tilde C_1^2) = - 2$,
$\overline {\text{NE}}(W)$ is generated by $C_1$ and an extremal ray
$\tilde \ell$ with respect to $K_W$.
We may assume that the strict transform $\ell$ of $\tilde \ell$ on $V$
is a $(- 1)$-curve,
because $(K_V \cdot \ell ) \le (K_W \cdot \tilde \ell ) < 0$ and we may take
an $\tilde \ell $ which passes through a singular point of $W$ if
$(\tilde \ell ^2) = 0$. We set $\bar \ell = f_*(\ell)$.

There are 2 possibilities: $\ell \ne E$ or $\ell = E$.
We assume first that $\ell \ne E$.
Since $(\bar \ell ^2) > 0$ and $(\tilde \ell ^2) \le 0$,
we have $(C_1 \cdot \ell ) > 0$. But $((K_V + D^{\#}) \cdot \ell) < 0$,
hence $(C_1 \cdot \ell ) = 1$ and $(E \cdot \ell) = 0$.
In [MT, Lemma 3.4] it is proved that $\ell$ does not intersect
$F$ except $C_1$.
Then by [MT, 4.1], $\ell $ intersects an irreducible component $D_1$ of
$D - F$.  If $(D_1^2) = - 2$, then $(D_1 \cdot \ell) = 1$,
because $(\tilde \ell ^2) \le 0$, hence
$\vert C_1 + 2\ell + D_1 \vert$ gives a permissible pencil
of rational curves.
If $(D_1^2) = - 3$, then we derive a contradiction with
$(\bar \ell ^2) > 0$ by the argument in [MT, p.295, up to line 17].

We shall assume that $\ell = E$ in the following.
Then we have $(\tilde E^2) \le 0$ and $(E^2) = - 1$.
The calculation in [MT, Lemma 3.1] shows that the combination of
singularities of $\bar V$ along $\bar E$ are as in one of the cases
(1), (2), (3), (5), (7), (9), (10), (11) or (13).
In the cases (1) with $r = 1$, (2), (5) or (9),
if $D_i = S_i \cap S_j$ in the notation of Theorem 5.3, then $D_i$
is contracted in $S_j$ to a rational double point.
Then $Z_s$ is Gorenstein, and so is $Z$,
and the claim of the theorem is clear.

In the case (1) with $r > 1$,
we consider a linear pencil $\Lambda = \vert C_1 + 2E + C_2 \vert$.
Then $A_r = C_3$ becomes a double section, and
there is only one $(- 1)$-curve on each singular fiber by [MT, 5.3.1.(i)].
Let $L$ be a singular fiber containing $A_{r-1}$, the irreducible
component of $F_3$ next to $C_3$.
By the argument in [MT, p.300, lines 5-11],
since $(A_r \cdot A_{r-1}) = 1$, we have $(A_r \cdot \ell) = 0$,
and the multiplicity of $A_{r-1}$ in $L$ is 2,
hence the graph of $L$ is either the bottom one on [MT, p.299]
or that on [MT, p.300].
In the former case, $D - F$ has a fork except when $s = 1$.  But
if $s = 1$, then $(B_1^2) = - (r + 1)$ in the notation there,
and the three curves at the
right end cannot be contracted to a singular point of the type
described in Corollary 4.3.
In the latter case, similar argument holds for the $s + 1$ curves at
the right end.

In the remaining cases, we shall first prove that
there exists a $(- 1)$-curve $h$
which does not meet $F - C_1$ and
such that $(h \cdot C_1) = 1$.
This is done in the cases (3), (7) and (13)
at [MT, 5.5 and the latter half of 5.3.2],
and in the cases (10) and (11) at
[MT, 5.6 and 5.7 using 5.2 similarly as 5.4].

We shall find an irreducible component $D_1$ of $D$ such that
$(D_1^2) = - 2$ and $(h \cdot D_1) = 1$.  Then the pencil
$\vert C_1 + 2h + D_1 \vert$ is permissible.
Let $\tilde h = g_*(h)$. In [MT, Lemma 6.1] one calculates
that $(\tilde h^2) < 1$ in our cases.
By the first argument of [MT, 6.2], there exists an irreducible
component $D_1$ of $D$ which meets $h$ and such that
$(D_1^2) = - 2$ or $- 3$.  Let $\mu = (h \cdot D_1)$.  In the latter case,
we have $\mu = 1$.  In the former case,
if we write $g^*\tilde h = h + \alpha _1D_1 + ...$ (instead of the equation
at [MT, p.306, line 10]), then
$0 = (g^*\tilde h \cdot D_1) \ge \mu - 2\alpha _1$, and
$(\tilde h^2) \ge - 1 + \mu ^2/2$, hence $\mu = 1$.
Then the rest is the same as in the last part of the proof in the case
in which $\ell \ne E$.  Q.E.D.
\enddemo

\proclaim{5.5. Flip Theorem}
Let $\varphi : X \to Z$ and $g : Z \to \Delta$ be morphisms which
satisfy the conditions in (5.1).
Then there exists a morphism $\varphi ^+: X^+ \to Z$
which satisfies the remaining conditions of (5.1), i.e.
the flip of $\varphi$.
Moreover, if $\text{Exc}(\varphi)$ is a geometrically irreducible curve,
then so is $\text{Exc}(\varphi ^+)$
\endproclaim

\demo{Proof}
We may assume that $A$ is a complete discrete valuation ring whose residue
field is algebraically closed. In fact, the existence of the flip
is equivalent to the surjectivity of the homomorphisms
$$
\varphi _*\Cal O_X(m_0K_{X/\Delta}) \otimes
\varphi _*\Cal O_X(nm_0K_{X/\Delta})
\to \varphi _*\Cal O_X((n + 1)m_0K_{X/\Delta})
$$
for a positive integer $m_0$ and all positive integers $n$,
and is stable by the flat base change.

Let  $r(\varphi)$ be the maximum of the indices of
$K_{X/\Delta}$ at the points on $\text{Exc}(\varphi)$.
We shall prove the theorem by induction on $r(\varphi)$.  Let $r_0$ be a
positive integer, and assume that
the theorem is true for any $\varphi$ such that $r(\varphi) < r_0$.
Let us fix a $\varphi$ such that $r(\varphi) = r_0$.

Let $Q \in Z$ be a point in the image of $\text{Exc}(\varphi )$. Since the
existence of the flip is a local question, we replace $Z$ by its
completion at $Q$, and then $X$, etc.,
by its fiber product over $Z$. We note
that the properties (1.1) are preserved under this replacement.

Let $e$ be the number of irreducible components of
$\text{Exc}(\varphi)$.  Then we have $\rho (X/Z) = e$
as in the proof of Theorem 5.3.
We may assume that $e = 1$ in the following.
In fact, if we have our theorem for $e = 1$, then applying
the minimal model program to $\varphi : X \to Z$,
we obtain a sequence of flips whose exceptional loci are
irreducible curves, which terminates by Theorem 5.2.  So we obtain
$\varphi ^{\#}:X^{\#} \to Z$ such that $K_{X^{\#}/\Delta}$ is
$\varphi ^{\#}$-nef.
Let $\psi _s: X^{\#}_s \to X^+_s$ be the contraction morphism
of all the curves whose intersection numbers with $K_{X^{\#}/\Delta}$ are 0.
Then we have $R^1\psi _{s*}\Cal O_{X^{\#}_s} = 0$, and $\psi _s$
can be extended to a morphism $\psi : X^{\#} \to X^+$ as in the proof of the
contraction theorem.  Since we can prove that
$R^1\psi _{s*}\Cal O_{X^{\#}_s}(mK_{X^{\#}/\Delta}) = 0$ for any integer $m$,
$X^+$ also satisfies
the conditions in (1.1).  Therefore, the induced morphism
$\varphi ^+:X^+ \to Z$ gives the flip of $\varphi$.

Let $C = \text{Exc}(\varphi)$, and
$\Lambda (\varphi)= \{P \in C: \text{ the index of
$K_{X/\Delta}$ at $P$ is $r_0$}\}$.
Let $n(P)$ be the axial multiplicity at $P \in \Lambda (\varphi)$,
and let $n(\varphi) = \sum _{P \in \Lambda (\varphi)}n(P)$
be the {\it total axial multiplicity}.

First, we shall prove our theorem under the
additional assumption that all the points $P \in \Lambda (\varphi)$ are of
simple type, and we proceed by induction on $n(\varphi)$.
We note that any non-normal point of $X_s$ is of simple type by
Theorem 4.1.

\proclaim{5.6. Lemma}
Let $P$ be a point of index $r_0$ on $C$, $\mu : X' \to X$
the standard blowing up at $P$, and $C'$  the strict transform of $C$.
Then $(K_{X'/\Delta}\cdot C') \le 0$.
\endproclaim

\demo{Proof}
There are three cases:
\roster
\item"(1)" $C \subset S_1 \cap S_2$ for two irreducible components $S_1$
and $S_2$ of $X_s$,

\item"(2)" $C \subset S_1$ for the only one irreducible component
$S_1$ of $X_s$ and $C \cap D_1 \ne \emptyset$,

\item"(3)" $C \subset S_1$ for the only one irreducible component
$S_1$ of $X_s$ and $C \cap D_1 = \emptyset$.
\endroster
Let $E$ be the exceptional divisor of $\mu$, and
$S'_i$ the strict transforms of the $S_i$ by $\mu$.

Case (1). Since $C$ is contractible, we have
$(S_i\cdot C) < 0$
for $i = 1, 2$. Since $(X_s\cdot C) = 0$,
there exists another irreducible component $S_3$ such that $(S_3\cdot C) > 0$.
Since the index of $K_{X/\Delta}$ at the point
$S_1 \cap S_2 \cap S_3$ is 1, $P$ is not contained in $S_3$.
Since the contraction of $C$ in $S_1$ yields a log terminal singularity,
$P$ is the only singular point of $S_1$ on $C$.
Then $C'$ intersects
two irreducible components $E$ and $S'_3$ of $X'_s$
besides $S'_1$ and $S'_2$.
Hence we have  $(S'_i\cdot C') = - 1$ for $i = 1, 2$, and
$(K_{X'/\Delta}\cdot C') = 0$, i.e., $C'$ is a $(- 1, - 1)$-curve.

Case (2). Since $C$ contracts to a log terminal singularity
inside $S_1$, it follows that $C \cap D_1$
consists of a single point $P_1$, and
there are at most another singular point $P_2$ of $S_1$ on $C$;
if $S_1$ is nonsingular along $C \setminus \{P_1\}$,
we let $P_2$ be an arbitrary point of $C \setminus \{P_1\}$.
Let $r_i \ge 1$ be the indices of $K_{X/\Delta}$ at the points $P_i$
for $i = 1, 2$.
Let $S_2$ be the other irreducible component of $X_s$ which intersects $C$,
and $D$ the irreducible component of $D_1$ which contains $P_1$.
We set $\frac{1}{r_1}(1, a_1)$ the type of the singularity of $S_1$
at $P_1$.
Let  $\alpha :S''_1 \to S_1$  be the minimal resolution, and let
$D''$ and $C''$ be the strict transforms of
$D$ and $C$ by $\alpha $, respectively.
Then $C''$ is a $(- 1)$-curve, and the exceptional divisors
$F_{i,j}$ ($i=1,2$ and $1 \le j \le d_i$ for some $d_i \ge 0$)
of $\alpha $ over $P_i$ satisfy
$$\align
(D''\cdot F_{1,1}) &= (F_{1,1}\cdot F_{1,2})
= ... = (F_{1,d_1-1}\cdot F_{1,d_1})= 1 \\
(C''\cdot F_{2,1}) &= (F_{2,1}\cdot F_{2,2})
= ... = (F_{2,d_2-1}\cdot F_{2,d_2})= 1
\endalign$$
under the suitable numbering.
There are two subcases:
\roster
\item"(2a)"  $P_2$ is a rational double point of type A and
$(C''\cdot F_{1,j_0}) = 1$ for some $j_0$,
and all other intersections are 0,

\item"(2b)" $(C''\cdot F_{1,d_1}) = 1$,
and all other intersections are 0.
\endroster

In the subcase (2a), since $r_2 = 1$, we have $r_0 = r_1$ and $P = P_1$.
Then $E_1 = E \cap S'_1$  is the image of  $F_{1,1}$, and
it is clear that $(K_{X'/\Delta}\cdot C') \le 0$, where the equality holds
if and only if  $j_0 = 1$.

In the subcase (2b), let $\frac{1}{h_2r_2^2}(b_2h_2r_2 - 1, 1)$
be the type of $P_2 \in S_1$ for some positive integers $b_2$ and
$h_2$ such that $b_2 < r_2$ (cf. Corollary 4.3).
We also set $a_1b_1 = c_1r_1 + 1$ for some integers
$b_1$ and $c_1$ such that $0 < b_1 \le r_1$.
The singularity of
$S_1$ at $P_1$ (resp. $P_2$) can be described as a toroidal embedding whose
lattice of 1-parameter subgroups is
$$
N = \Bbb Z^2 + \frac{1}{r_1}(1, a_1)\Bbb Z
\quad \text{ (resp. $\Bbb Z^2 + \frac{1}{h_2r_2^2}
(b_2h_2r_2 - 1, 1)\Bbb Z$)},
$$
and $D''$, $F_{1,d_1}$, $C''$
(resp. $C''$, $F_{2,1}$) correspond to
$(0, 1)$, $(b_1/r_1, 1/r_1)$, $(1, 0)$
(resp. $(1, 0)$, $((b_2h_2r_2 - 1)/h_2r_2^2, 1/h_2r_2^2)$ )
in  $N_{\Bbb Q}$. Then the total
transform of $C$ on $S''_1$ has coefficient
$b_1/r_1$ (resp. $(b_2h_2r_2 - 1)/h_2r_2^2$) at
$F_{1,d_1}$ (resp. $F_{2,1}$) and
$$
(C^2) = - 1 + b_1/r_1 + (b_2h_2r_2 - 1)/h_2r_2^2 < 0,
$$
because $C$ is contractible. The discrepancy coefficient of $F_{1,d_1}$
(resp. $F_{2,1}$) with respect to $K_{S_1} + D$ is given by
$$
(b_1+1)/r_1 - 1 - 1/r_1 = b_1/r_1 - 1
\text{ (resp. $b_2/r_2 - 1$)},
$$
hence
$$
((K_{S_1} + D)\cdot C) = - 1 - (b_1/r_1 - 1) - (b_2/r_2 - 1)
= 1 - b_1/r_1 - b_2/r_2 < 0.
$$
Thus
$$
b_1/r_1 + b_2/r_2 - 1/h_2r_2^2 < 1 < b_1/r_1 + b_2/r_2 = 1 + e/r_1r_2
$$
for a positive integer $e$. Hence
$1/h_2r_2^2 > 1/r_1r_2$, i.e., $r_1 > h_2r_2$. Therefore, $r_0 = r_1$ and
$P = P_1$.

The exceptional divisor $E_1 = E \cap S'_1$
is the image of $F_{1,1}$, which corresponds to $(1/r_1, a_1/r_1)$.
The new singularity at the
intersection $E_1 \cap C'$ is of type $\frac{1}{a_1}(1, - r_1)
= \frac{1}{a_1}(c_1, 1)$, and
$$
(K_{X'/\Delta}\cdot C') = ((K_{S'_1} + D' + E_1)\cdot C') =
- 1 - (c_1/a_1 - 1) - (b_2/r_2 - 1)
= 1/a_1r_1 - e/r_1r_2.
$$
Since $b_1r_2 = (r_2 - b_2)r_1 + e \equiv e (\text{mod } r_1)$,
we have $r_2 \equiv a_1e (\text{mod } r_1)$,
hence $r_2 \le a_1e$. Therefore, $(K_{X'/\Delta}\cdot C') \le 0$.

Case (3). If there is only one singular
point of index greater than one on $C$, then it is clear that
$(K_{X'/\Delta}\cdot C') \le 0$.
We assume that there are more than one such singularities in the
following. Then there are exactly 2 singular points, say $P_1$ and $P_2$,
because $C \subset S_1$ contracts to a rational singularity.
Let $\alpha :S''_1 \to S_1,C''$, and
$F_{i,j}$ ($i=1,2$ and $1 \le j \le d_i$ for some $d_i \ge 0$)
be as in the case (2) except the intersection number with
$D''$. As before, let
$\frac{1}{h_ir_i^2}(b_ih_ir_i - 1, 1)$
be the types of the $P_i \in S_1$ for some positive integers $r_i, b_i$ and
$h_i$ such that $b_i < r_i$ for $i = 1, 2$.
We set also $a_ib_i = c_ir_i + 1$ for some integers
$a_i$ and $c_i$ such that $0 < a_i < r_i$.
There are 2 subcases:
\roster
\item"(3a)"  $(C''\cdot F_{1,d_1}) = 1$,
and all other intersections are 0,

\item"(3b)"  $(C''\cdot F_{1,j_0}) = 1$ for some $j_0$
with $1<j_0<d_1$,
and all other intersections are 0.
\endroster

In the subcase (3a), $C$ is contracted to a log terminal singularity of type
A. We may assume $r_0 = r_1 \ge r_2, h_1 = 1$ and $P = P_1$.
Then
$$\align
(C^2) &= - 1 + (b_1r_1 - 1)/r_1^2 + (b_2h_2r_2 - 1)/h_2r_2^2 < 0 \\
(K_{X/\Delta}\cdot C) &= - 1 - (b_1/r_1 - 1) - (b_2/r_2 - 1) < 0,
\endalign$$
thus
$$
b_1/r_1 + b_2/r_2 - 1/r_1^2 - 1/h_2r_2^2 < 1 < b_1/r_1 + b_2/r_2 =
1 + e/r_1r_2
$$
for a positive integer $e$. The exceptional divisor $E_1 = E \cap S'_1$
corresponds to
$((r_1 - a_1)/r_1^2, a_1/r_1^2) \in N_{\Bbb Q}$ as before. We have
$$
((K_{S'_1} + E_1)\cdot C') = - 1 - (c_1/a_1 - 1) - (b_2/r_2 - 1)
= 1/a_1r_1 - e/r_1r_2.
$$

Since $b_1r_2 = (r_2 - b_2)r_1 + e \equiv e (\text{mod } r_1)$,
we have $r_2 \equiv a_1e (\text{mod } r_1)$,
hence $r_2 \le a_1e$. Therefore, $(K_{X'/\Delta}\cdot C') \le 0$.

We shall prove that the subcase (3b) does not occur.
If $C$ is contracted to a log terminal singularity of type D,
then we have one of the following:
\roster
\item"(3b1)"  $d_1 = 3$, $j_0 = 2$, $(F_{1,1}^2) = (F_{1,3}^2) = - 2$, and
$(F_{1,2}^2) = \ell$ for an integer $\ell \ge 3$,

\item"(3b2)"  $d_2 = 1$, $j_0 = 2$, $(F_{1,1}^2) = - 2$, and
$(F_{2,1}^2) = - 3$.
\endroster
But then, $P_1$ (resp. $P_2$) in the case (3b1) (resp. (3b2) )
is of type $\frac{1}{4\ell - 4}(2\ell - 1, 1)$
(resp. $\frac{1}{3}(1, 1)$), and
cannot be of type
$\frac{1}{hr^2}(bhr - 1, 1)$ for some positive integers $r, b$ and $h$
with $0 < b < r$.

Next, assume that a type E singularity appears after the contraction.
As for $P_2$, a similar consideration
as above shows that
$d_2 = 1$ and $(F_{2,1}^2) = - 4$, or
$d_2 = 2$ and $(F_{2,1}^2) = (F_{2,2}^2) = - 3$.
Then for $P_1$, we have  $j_0 = 2$, and one of the following holds
for an integer $\ell \ge 2$:
\roster
\item"(3b3)"  $d_1 = 3$, $(F_{1,1}^2) = - 2$, $(F_{1,2}^2) = - \ell$,
and $(F_{1,3}^2) = - 3$,

\item"(3b4)"  $d_1 = 3$, $(F_{1,1}^2) = - 2$, $(F_{1,2}^2) = - \ell$,
and $(F_{1,3}^2) = - 4$,

\item"(3b5)"  $d_1 = 3$, $(F_{1,1}^2) = - 2$, $(F_{1,2}^2) = - \ell$,
and $(F_{1,3}^2) = - 5$,

\item"(3b6)"  $d_1 = 4$, $(F_{1,2}^2) = - \ell$,
and $(F_{1,1}^2) = (F_{1,3}^2) = (F_{1,4}^2) = - 2$,

\item"(3b7)"  $d_1 = 4$, $(F_{1,1}^2) = (F_{1,3}^2) = - 2$, $(F_{1,2}^2) =
- \ell$, and $(F_{1,4}^2) = - 3$,

\item"(3b8)"  $d_1 = 4$, $(F_{1,1}^2) = (F_{1,4}^2) = - 2$,
$(F_{1,2}^2) = - \ell$, and $(F_{1,3}^2) = - 3$,

\item"(3b9)"  $d_1 = 5$, $(F_{1,2}^2) = - \ell$,
and $(F_{1,1}^2) = (F_{1,3}^2) = (F_{1,4}^2) = (F_{1,5}^2) = - 2$,

\item"(3b10)"  $d_1 = 6$, $(F_{1,2}^2) = - \ell$,
and $(F_{1,1}^2) = (F_{1,3}^2) =
(F_{1,4}^2) = (F_{1,5}^2) = (F_{1,6}^2) =- 2$.
\endroster
Then the only possible cases are (3b3) with $\ell = 5$, (3b4) with
$\ell = 3$, and (3b7) with $\ell = 6$.
But in these cases, we have
$$\align
(K_{X/\Delta}\cdot C) &= - 1 + 1/2 + 4/5 > 0 \\
(K_{X/\Delta}\cdot C) &= - 1 + 1/2 + 2/3 > 0 \\
(K_{X/\Delta}\cdot C) &= - 1 + 1/2 + 6/7 > 0,
\endalign$$
a contradiction. Q.E.D.
\enddemo

{\it Proof of Theorem 5.5 continued.}
If $(K_{X'/\Delta}\cdot C') < 0$, we put $X^{(0)} = X'$.
Otherwise, we have $(K_{X'/\Delta}\cdot C') = 0$.
Then there exists a projective birational morphism
$\psi : X' \to Z'$ over $\Delta$ which contracts only $C'$
as in the contraction theorem.
We shall construct a {\it flop} of $\psi$:
$$
X' @>\psi >> Z' @<\psi ^+<< X^{(0)}.
$$
In the case (1) of Lemma 5.6, the blowing up
of $C'$ and the contraction in the other direction give us the flop.
In other cases, since  $R^1\psi _*\Cal O_{X'}(mK_{X'/\Delta}) = 0$
for any integer $m$, $Z'$ satisfies the conditions in (1.1) except
that the irreducible components of $Z'_s$
are not necessarily $\Bbb Q$-Cartier.
Let $Q' = \psi (C')$ and $r'$ the index of $K_{Z'/\Delta}$ at $Q'$.
As in Theorems 4.1 or 4.4, we can prove that
the completion of the canonical cover of $\Cal O_{Z',Q'}$
with the $\mu _{r'}$-action is isomorphic to
that of $A[x, y, z]/(F')$ with $F' = xy + \tau G'(z^{r'})$
for some semi-invariant coordinates $x, y, z$ and a generator $\tau$ of
the maximal ideal of A, where the action is given by
$x \mapsto \zeta ^{a'} \otimes x$, $y \mapsto \zeta ^{-a'} \otimes y$, and
$z \mapsto \zeta \otimes z$
for some positive integer $a'$ such that $(r', a') = 1$.
Since the strict
transforms of $S_1$ and $E$ on $Z'$ are not $\Bbb Q$-Cartier at $Q'$,
$G'$ is not invertible. So if $\psi $ is obtained from the blowing up
at the ideal $(x, \tau)$, then the blowing up at $(y, \tau)$
gives us the flop $\psi ^+: X^{(0)} \to Z'$.

Now we apply the minimal model program to the induced morphism
$f^{(0)}: X^{(0)} \to Z$.  We note that $\rho(X^{(0)}/Z) = 2$.
Since $K_{X^{(0)}/\Delta}$ is not relatively nef over $Z$,
there exists a contraction morphism $\varphi ^{(0)}: X^{(0)} \to Z^{(0)}$
with respect to an extremal
ray of $\overline{\text{NE}}(X^{(0)}/Z)$.
In the case $(K_{X'/\Delta} \cdot C') < 0$,
we take $\varphi ^{(0)}$ to be the contraction of $C'$.

If $\varphi ^{(0)}$ is a small contraction,
then the flip $\varphi ^{(0)+}: X^{(1)} \to Z^{(0)}$ exists, since
$r(\varphi ^{(0)}) \le r_0$, and if the equality holds, then
$n(\varphi ^{(0)}) < n(\varphi)$.  Let $\varphi ^{(1)}: X^{(1)} \to Z^{(1)}$
be the contraction morphism with respect to an extremal
ray of $\overline{\text{NE}}(X^{(1)}/Z)$.
If it is a small contraction again, then the flip
$\varphi ^{(1)+}: X^{(2)} \to Z^{(1)}$
exists by the same reason. We continue this process. By Theorem 5.2, we
obtain a divisorial contraction
$\varphi ^{(s)}: X^{(s)} \to Z^{(s)} = X^+$ for a nonnegative
integer $s$, where $s = 0$ may occur only if $(K_{X'/\Delta} \cdot C')
= 0$.
By Theorem 5.3, $X^+$ satisfies the conditions in (1.1).

We shall prove that the composite birational map $X \dasharrow X^+$ is not
the identity. Then the induced morphism $\varphi ^+: X^+ \to Z$ gives us the
flip of $\varphi $, because $\rho (X^+/Z) = 1$ and
$K_{X^+/\Delta}$ must be $\varphi ^+$-ample.
Let $v$ be a discrete valuation of the fraction field of $X$
which has the center $C'$ on $X'$. Let $d$ and $d^+$ be the discrepancy
coefficients of $v$ over $X$ and $X^+$, respectively.
If $s = 0$, then $d < d^+$,
because $C'$ is not contained in $E$, while the flopped curve is
contained in the strict transform of $E$.
If $s > 1$, then $d < d^+$,
because a flip or a divisorial contraction increases the discrepancy.
Therefore, the theorem is proved under the additional assumption that
all the singular points of index $r_0$ are of simple type.

Next, we consider the general case.  By the proof of Lemma 5.6, we have
to consider only the case (3) there.  Let $P \in \Lambda (\varphi)$.
By Theorem 4.1, the completion of $X$ at $P$ has an equation
$xy + G_P(z^{r_0}) = 0$ in the $\mu _{r_0}$-quotient of
the completion of $\text{Spec } A[x, y, z]$.
Since $X_{\eta }$ is smooth, the generic fiber of the 1-dimensional scheme
defined by
$G_P(w) = 0$ in the completion of $\text{Spec } A[w]$
is geometrically reduced.

Let $\tilde A$ be complete discrete valuation ring which is finite over $A$
and such that $\tilde \Delta = \text{Spec } \tilde A$
dominates all the irreducible components of these 1-dimensional schemes
for all $P \in \Lambda (\varphi)$.
Let $\tilde X = X \times _{\Delta} \tilde \Delta$,
$\tilde Z = Z \times _{\Delta} \tilde \Delta$,
$\sigma : \tilde X \to X$ and $\tilde \varphi : \tilde X \to \tilde Z$
the natural morphisms, $\tilde P = \sigma ^{- 1}(P)$ for
$P \in \Lambda(\varphi)$, and
$\tilde \tau$ be a generator of the maximal ideal of $\tilde \Delta$,
Since $X_s$ is normal,
the completion of $\tilde X$ at $\tilde P$ is described by an
equation $xy + \prod _j G_{P,j}(z^{r_0}) = 0$
where the $G_{P,j}(\text{mod }\tilde \tau \tilde A[z^{r_0}])
\in z^{r_0}k[z^{r_0}]$
contain linear terms in $z^{r_0}$.
By the argument in [K3, p.116], there exists a projective birational morphism
$\lambda :Y \to \tilde X$ which contracts only a finite number of curves and
such that
$K_{Y/\tilde \Delta} = \lambda ^*K_{\tilde X/\tilde \Delta}$,
the induced morphism
$g: Y \to \tilde \Delta$ satisfies the conditions in (1.1), and that
$Y$ has only singularities of simple type.
If we apply the minimal model program to the small morphism
$\tilde \varphi \circ \lambda: Y \to \tilde Z$,
then there exists a projective birational morphism
$\tilde \varphi ^+: Y^+ \to \tilde Z$ such that
$K_{Y^+/\tilde \Delta}$ is $\tilde \varphi ^+$-nef
by the first part of the proof.
Then there exists a positive integer $m_0$ such that the natural
homomorphisms
$$
\tilde \varphi _*\Cal O_{\tilde X}(m_0K_{\tilde X/\tilde \Delta}) \otimes
\tilde \varphi _*\Cal O_{\tilde X}(nm_0K_{\tilde X/\tilde \Delta})
\to \tilde \varphi _*\Cal O_{\tilde X}((n + 1)m_0K_{\tilde X/\tilde \Delta})
$$
are surjective for any positive integer $n$.
By the flat base change theorem, we deduce the existence of the flip
of $\varphi$.
Q.E.D.
\enddemo

As explained in [KMM], the cone, contraction, and flip theorems combined
yield the following:

\proclaim{5.7. Minimal Model Theorem}
Let $f^{(0)}: X^{(0)} \to \Delta $ be a semistable projective
morphism of relative dimension 2.
Then there exist a projective flat morphism
$f: X \to \Delta $ from a $\Bbb Q$-factorial scheme $X$
satisfying the conditions in \rom{(1.1)}
and a birational map $\alpha : X^{(0)} \dasharrow X$
over $\Delta $ such that one of the following holds:
\roster
\item"(a)" $K_{X/\Delta}$ is $f$-nef,

\item"(b)" $\rho (X/\Delta ) = 2$,
and there exist a projective flat morphism
$g: Z \to \Delta $ of relative dimension 1
and a projective surjective morphism
$\varphi : X \to Z$ over $\Delta $ such that
$- K_{X/\Delta}$ is $\varphi $-ample,

\item"(c)" $\rho (X/\Delta ) = 1$, and $- K_{X/\Delta}$ is $f$-ample.
\endroster
\endproclaim

In the case (a), if $\kappa (X^{(0)}_{\eta}) = 0$, then
$\Cal O_X(12K_{X/\Delta}) \in f^*\text{Pic}(\Delta)$
by the upper semicontinuity theorem.

\proclaim{5.8. Theorem}
In the case (a) of Theorem 5.7, if $\kappa (X^{(0)}_{\eta}) = 2$, then
there exists a positive integer $m_0$
such that the natural homomorphism
$$
f^*f_*\Cal O_X(m_0K_{X/\Delta}) \to \Cal O_X(m_0K_{X/\Delta})
$$
is surjective.  In particular, the relative canonical ring
$$
\bigoplus _{m=0}^{\infty} f^{(0)}_*\Cal O_{X^{(0)}}(mK_{X^{(0)}/\Delta})
$$
is finitely generated over $A$.
\endproclaim

\demo{Proof}
We may assume that $A$ is a complete discrete valuation ring whose residue
field is algebraically closed.
First, we shall prove that $\vert m_1K_{X_s} \vert$ is free for a positive
integer $m_1$.
The proof is similar to that in [K5, middle of p.356 to 357].
But we do not use the log canonical cover since it is not necessarily normal.
Let $\nu _i = \kappa (S_i, K_{S_i} + D_i)$
for the irreducible components $S_i$ of $X_s$.
We may assume that $D_i \ne 0$ for all $i$ when we glue sections of
pluricanonical sheaves.
Under this assumption, if $\nu _i = 0$, then $S_i$ is birationally
equivalent to a ruled surface over a curve $C_i$ of genus $g_i$.
If $g_i > 0$, then
$D_i$ contains either 2 sections or a double section which may be an
inseparable section. Hence $g_i \le 1$.

We claim that $D_i$ is connected if $\nu _i = 0$ and $S_i$ is rational.
If we apply the log minimal model program to the pair $(S_i, 0)$,
we obtain a projective birational morphism $\mu _i: S_i \to S'_i$ to a
surface $S'_i$ with only log terminal singularities
such that $\rho (S'_i) = 1$, or that $\rho (S'_i) = 2$ and
there is a surjective morphism $\varphi _i: S'_i \to C_i$ whose
fibers are irreducible.
Then $D'_i = \mu _{i*}(D_i)$ is connected.
Tracing back the morphism $\mu _i$, we conclude that $D_i$ is also
connected.

If $\nu _i = 0$, then $((K_{S_i} + D_i) \cdot \Gamma) = 0$ for any
$\Gamma \subset D_i$, and there are the following possibilities:
\roster
\item"(1)" $S_i$ is a rational surface, and $D_i$ is a nonsingular elliptic
curve,

\item"(2)" $S_i$ is a rational surface, and $D_i$ is a rational curve
with a node or a cycle of nonsingular rational curves,

\item"(3)" $S_i$ is a rational surface, and $D_i$ is a nonsingular
rational curve on which there are 4 ordinary double points of $S_i$,
or a rod of nonsingular rational curves each of whose 2 end components
carries 2 ordinary double points of $S_i$,

\item"(4)" $S_i$ is a rational surface, and $D_i$ is a nonsingular rational
curve on which there are 3 singular points of $S_i$,

\item"(5)" $S_i$ is birationally equivalent to an elliptic ruled surface,
and $D_i$ consists of 2 disjoint nonsingular elliptic curves which are
sections of the ruling,

\item"(6)" $S_i$ is birationally equivalent to an elliptic ruled surface,
and $D_i$ is a nonsingular elliptic curve which is a double section of the
ruling.
\endroster
In the case (4), the log indices of $K_{S_i} + D_i$ at
the singular points are $(2, 3, 6)$,
$(2, 4, 4)$, or $(3, 3, 3)$.
Therefore as in [K5], we obtain canonical sections in
$H^0(D_i, 12(K_{S_i} + D_i))$, hence a section in
$H^0(B(0), 12m_2K_{X_s})$ for a positive integer $m_2$, where $B(0)$ is the
union of the $S_i$ with $\nu _i = 0$. Then the rest is the same as in [K5]
and $\Cal O_{X_s}(m_1K_{X_s})$ is generated by global sections for a
positive integer $m_1$.

Next, we shall prove that $H^1(X_s, m_0K_{X_s}) = 0$ for a positive multiple
$m_0$ of $m_1$.  Then the sections in  $H^0(X_s, m_0K_{X_s})$ are extended
to $X$, and the proof of the theorem is completed.
The free linear system $\vert m_1K_{X_s} \vert$ gives a morphism
$\psi : X_s \to Y_s$ with connected fibers
to a 2-dimensional projective scheme $Y_s$
such that $m_1K_{X_s} = \psi ^*H$ for an ample divisor
$H$ on $Y_s$.  Thus it is enough to prove that
$R^1\psi_*\Cal O_{X_s} = 0$.
The proof is similar to that of Theorem 2.3.
If $\nu _i = 2$, then $\psi _i = \psi \vert _{S_i}$ contracts only a
finite number of curves, and we have $R^1\psi _{i*}\Cal O_{S_i}(- D_i) = 0$.
For $\nu _i = 1$, we can use Lemma 2.1.
So we have to prove that $H^1(B(0), \Cal O_{B(0)}) = 0$.

Let $K^0 = \oplus _i \Cal O_{S_i},
K^1 =  \oplus _{i > j} \Cal O_{D_{ij}}$, and
$K^2 = \oplus _{i > j > k} \Cal O_{P_{ijk}}$,
where $D_{ij} = S_i \cap S_j$,
$P_{ijk} = S_i \cap S_j \cap S_k$, and
the summations are taken for all $i, j, k$ such that
$\nu _i = \nu _j = \nu _k = 0$.
Then we have an exact sequence
$$
0 \to \Cal O_{B(0)} \to K^0 \to K^1 \to K^2 \to 0.
$$
Since there are components $S_i$ such that $\nu _i = 2$,
we can check that the homomorphism
$\oplus _i H^1(S_i, \Cal O_{S_i}) \to
\oplus _{i > j} H^1(D_{ij}, \Cal O_{D_{ij}})$ is injective, and that the
dual graph of intersections of irreducible components of $B(0)$ has no
cycle.  Then by the spectral sequence
$E_1^{p, q} = H^q(K^p) \Rightarrow H^{p+q}(B(0), \Cal O_{B(0)})$,
we obtain $H^1(B(0), \Cal O_{B(0)}) = 0$.
Q.E.D.
\enddemo

\enddocument